\documentclass[%
 reprint,
nofootinbib,
 amsmath,amssymb,
 aps,
prd,
]{revtex4-2}
\makeatletter
\count\footins=1000
\makeatother
\usepackage[T1]{fontenc}
\usepackage[utf8]{inputenc}
\usepackage{lmodern}
\usepackage{graphicx}
\usepackage{subfigure}
\usepackage{dcolumn}
\usepackage{bm}
\usepackage{hyperref}

\begin{document}

\preprint{APS/123-QED}

\title{Double twist knots and lattice paths}

\author{Aditya Dwivedi}
\email{aditya.dwivedi@iitb.ac.in}
\author{Ramadevi Pichai}%
 \email{ramadevi@iitb.ac.in}
\affiliation{Department of Physics, Indian Institute of Technology Bombay,\break Powai, Mumbai 400076, India}



\begin{abstract}
In this work, we explore the combinatorics arising from the quiver generating series of the unreduced $r$-colored HOMFLY-PT polynomial $\bar{P}_r(a,q)$ for some twist-knots and  double twist knots. By taking the limit $a = 0$ and $q = 1$, we indeed obtain lattice path models for these knots.
\end{abstract}

\maketitle
\tableofcontents

\section{Introduction}

 Knot invariants such as the HOMFLY-PT polynomial play a central role in low-dimensional topology and reveal deep connections with combinatorics, algebra, and geometry. From the point of view of both mathematics and physics, knot invariants are expected to have a counting interpretation. On the physics side, they describe the counting of BPS states \cite{kucharski2017bps}. Similarly, on the mathematical side, these invariants carry homological structure \cite{gorsky2010q}. In recent years, combinatorics has become an important tool in understanding this counting interpretation of knot invariants \cite{Kucharski:2016rlb}, and lattice path counting has played a major role in this progress. 
\par
These path counting models are the consequence of the  knot-quiver correspondence \cite{kucharski2017knots}. Interestingly, a lattice path counting was applied to the study of the $r$-colored HOMFLY-PT polynomial $P_r(a,q)$ of torus knots with symmetric color: $%
\begin{array}{|c|c|c|c|}
\hline
1 & 2 & \cdots & r \\
\hline
\end{array}%
$ \cite{Panfil:2018sis, Stosic:2024vbd, bizley1954derivation} . Particularly, it was shown that the quiver generating function of unreduced colored HOMFLY-PT polynomials $\bar{P}_r(a,q)$\footnote{We choose the normalization
$\bar{P}_r(a,q)(K)=P_r(a,q)\,\bar{P}_r(\mathrm{unknot})
=P_r(a,q)\,a^{-r}q^r\frac{(a^2;q^2)_r}{(q^2;q^2)_r}$.} for torus knots can be interpreted as a lattice path counting model in the positive quadrant. An important structural feature that appears is that these generating functions can be mapped to the generating function of a lattice path model under a line of rational slope passing through the origin. This idea was later extended to some non-torus knots, where the corresponding models involve a line of rational slope shifted from the origin \cite{Dordevic:2024jha}. The success of the lattice path model for such knots \cite{Panfil:2018sis, Stosic:2024vbd, Dordevic:2024jha} motivates us to extend the lattice model investigation for double twist knots belonging to a  $3$-pretzel knot family. 

\begin{figure}[t]
    \centering
    \includegraphics[width=1.05\columnwidth]{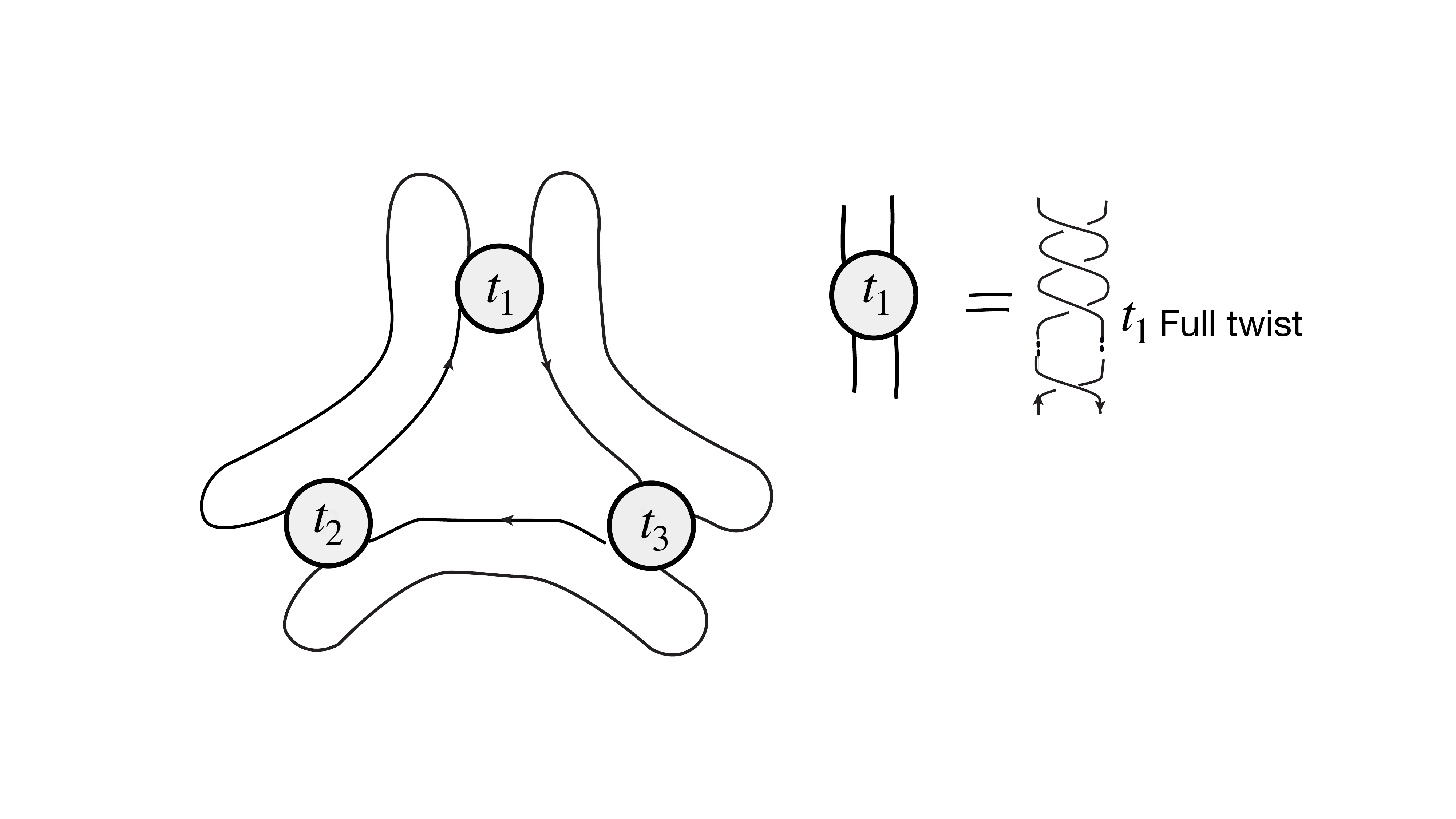}
    \caption{$L(t_1,t_2,t_3)$ knot.}
    \label{fig:knot}
\end{figure}

$3$-pretzel knot $L(t_1,t_2,t_3)$  has three twist regions where $t_1, t_2, t_3$ represent the number of full twist in each region (see Figure \ref{fig:knot}). 
This pretzel knot family contains twist knots $K_p$  as well as double twist knots $K(p_1,p_2)$ as subfamilies. More explicitly, the twist knot $K_{p<0} \equiv L(2p+1,1,1)$ in the $3$-pretzel knot notation. Similarly, the 3-pretzel form for the positive full-twist knots is $K_{p>0} \equiv L(2p-1,1,1)$. 
Further, the double twist knots in the $3$-pretzel notation are represented as $K(p_1+1,p_2+1)\equiv L(2p_1+1,2p_2+1,1)$. Note that the twist knots 
$K_{p>0}\equiv K(p,1)$ is the notation of double twist knots.
Our aim is to attempt lattice path counting for the twist knots and double-twist knots as a function of twisting parameter(s).

\par
Lattice path counting is dependent on the knowledge of quiver generating series for the knots. Such quiver generating series are known for torus knots, twist knots, and some double twist knots \cite{kucharski2017knots, Singh:2023qpw, Banerjee:2020dqq, Chauhan:2025iwt}. The lattice path results for torus knots $T(2,2p+1)$ dependent on twist parameter $p$, suggest that we should find the lattice paths for twist knots and double twist knots dependent on the full-twist parameters.

Interestingly, the recent paper by one of the authors \cite{Chauhan:2025iwt} on generating the quiver presentation for 3-pretzel knot family gives a tangible idea to systematize the lattice counting.  
In such a set-up, we first need to obtain the augmented quiver $\tilde{Q}^+_{K_1}$ for the knot $K_1$ with lower number of twists. For instance, $4_1 \equiv K_{p=-1}$ is the lowest knot in the twist knot family. 
Then, invoking {\it twisting on knots} to the {\it unlinking/linking} operations on the augmented quiver gives the quiver generating series for the tower of knots $\{K_i\}$'s. This is useful to deduce the lattice counting dependence on the twist parameter(s). We succeeded in deducing the lattice path interpretation for some twist knots and double twist knots with different framing $f$.

This paper is organized as follows: In Section \ref{Sec:Math_prelm}, we will first briefly review the salient aspects of the knot-quiver correspondence. We then present a choice of specialization to convert multi-variable quiver generating series to one-variable series, which will allow a lattice path counting interpretation. We also give a concise summary of the unlinking/linking operation on the quivers, which is essential to obtain the quiver adjacency matrix for the tower of knots obtained by full-twist operations.
In Section \ref{Sec:Examples}, we briefly review the lattice path counting results for torus knots  $T(2,2p+1)$ \cite{Panfil:2018sis, Stosic:2024vbd}. We then present our new results for some twist knots and double twist knots which belong to the  $3-$pretzel knot $L(t_1,t_2,t_3)$ family. We have put forth the lattice path counting results for these knots with different framing. From these results, we observe the lattice path counting for twist knots as a function of the twisting parameter. The concluding Section \ref{Sec:conclusion}  gives the summary of our results and possible future direction.
\section{Mathematical Preliminaries}\label{Sec:Math_prelm}
We begin with the relevant mathematical details essential to obtain lattice counting series. We now present a brief overview of the knot-quiver correspondence, its lattice path interpretation, and augmented quivers of knots based on \cite{kucharski2017knots, Stosic:2024vbd, Chauhan:2025iwt}. 
\subsection{The knot-quiver correspondence}
The knot-quiver correspondence basically connects the $P_r(a,q)$ generating series with quiver data \cite{kucharski2017knots}. For a knot $K$, let $Q_K (C,\mathbf{x})$ denote the associated quiver. Here, $C$ is the symmetric adjacency matrix of dimension $k\times k$ and $\mathbf{x}$ are the generators of the quiver series. In the reduced setting, the generating series of colored HOMFLY-PT polynomials can be written in quiver form using $C$ together with the grading vectors
\begin{equation}    
\mathbf{a}=(a_1,\ldots,a_k), \, \mathbf{q}=(q_1,\ldots,q_k), \label{aqgradreduced}
\end{equation}
which record the $(a,q)$-gradings of the colored HOMFLY-PT generators. 
For reduced HOMFLY-PT polynomials, the correspondence takes the form
\begin{equation}
\sum_{r\ge 0}\frac{P_r(a,q)(K)}{(q^2;q^2)_r}\,x^r
=
P_C(x_1,\dots,x_k)\Big|_{x_i=a^{a_i}q^{\,q_i-C_{ii}}x},
\tag{9}
\end{equation}
where the quiver generating series is
\begin{equation}
P_C(x_1,\dots,x_k)
=
\sum_{d_1,\dots,d_k}
\frac{
(-q)^{\sum_{i,j}C_{ij}d_i d_j}
x_1^{d_1}x_2^{d_2}\cdots x_k^{d_k}
}{
\prod_{i=1}^{k}(q^2;q^2)_{d_i}
}.
\end{equation}
The quiver matrix for mirror image of $K$ is described by replacing $C$ with $I_{k\times k}-C$.

The adjacency matrix $\bar C$ for the unreduced colored HOMLFY-PT $\bar{P}_r(a,q)$ is needed for the study of the lattice-path model. The explicit form of $\bar C$ \cite{Stosic:2024vbd} is
\begin{equation}
\bar{C}
=
C \otimes
\begin{pmatrix}
1 & 1 \\
1 & 1
\end{pmatrix}
+
\begin{pmatrix}
1 & 0 & 1 & 0 & 1 & 0 & \cdots & 1 & 0 \\
0 & 0 & 1 & 0 & 1 & 0 & \cdots & 1 & 0 \\
1 & 1 & 1 & 0 & 1 & 0 & \cdots & 1 & 0 \\
0 & 0 & 0 & 0 & 1 & 0 & \cdots & 1 & 0 \\
1 & 1 & 1 & 1 & 1 & 0 & \cdots & 1 & 0 \\
0 & 0 & 0 & 0 & 0 & 0 & \cdots & 1 & 0 \\
\vdots & \vdots & \vdots & \vdots & \vdots & \vdots & \ddots & \vdots & \vdots \\
1 & 1 & 1 & 1 & 1 & 1 & \cdots & 1 & 0 \\
0 & 0 & 0 & 0 & 0 & 0 & \cdots & 0 & 0
\end{pmatrix}.
\end{equation}
Note that the above matrix is of $2k \times 2k$ dimensions. Adding framing $f$ to the knot $K$ gives 
$\bar C_{K_f}$ whose matrix form in terms of $\bar C$ is 
$$\bar{C} \longmapsto \bar C_{K_f} =\bar{C} + f\,\mathbf{1}_{2k\times 2k}~.$$
Hence, the quiver generating series for unreduced colored HOMLFY-PT polynomials for such framed knots $K_f$ is as follows:
\begin{equation}
\sum_{r\geq 0}\bar{P}_r(K_f)x^r
= \bar P_{\bar C_{K_f}}(x_1,\ldots,x_{2k})\bigg|_{x_i=a^{\bar a_i}q^{\bar q_i-(\bar C_{K_f})_{ii}}x},
\label{kqc-unreduced-final}
\end{equation}
where
\begin{equation}
\begin{split}
\bar{P}_{\bar C_{K_f}}(x_1,\ldots,x_{2k})
=&\sum_{d_1,\ldots,d_{2k}}
(-q)^{\sum_{i,j}(\bar C_{K_f})_{ij}d_id_j} \\
&\times
\frac{x_1^{d_1}\cdots x_{2k}^{d_{2k}}}
{\prod_{i=1}^k (q^2;q^2)_{d_i}}.
\end{split}
\label{kqc-unreduced-series-final}
\end{equation}
Here, 
\begin{align}
\mathbf{\bar{a}}
&=(\bar a_1,\bar a_2,\ldots,\bar a_{2k}) \nonumber\\
&=(a_1+1,a_1-1,a_2+1,a_2-1,\dots,a_k+1,a_k-1), \nonumber\\[1mm]
\mathbf{\bar{q}}
&=(\bar q_1,\bar q_2,\ldots,\bar q_{2k}) \nonumber\\
&=(q_1+1,q_1+1,q_2+1,q_2+1,\dots,q_k+1,q_k+1)
\end{align}
are the degree vectors assigned to the vertices of the quiver 
$\tilde Q_{K_f}$ written in terms of the $(a,q)$
 gradings of the reduced quiver $Q_K$ (\ref{aqgradreduced}).

Our focus is to consider suitable one-variable specializations instead of the full multi-variable quiver series as discussed in Ref. \cite{Stosic:2024vbd}. Such a specialization will provide the generating functions from which we can possibly extract the lattice-path description for many knots.
Using the following shifting on the $\mathbf {\bar a}$:
\begin{equation}
a_{\max}=\max_{1\leq i\leq 2k} \bar a_i,
\qquad
\bar a_i'=a_{\max}-\bar a_i,
\qquad i=1,\dots,2k~,
\end{equation}
the following substitution was implemented for $T(2,2p+1)$ torus knots \cite{Stosic:2024vbd}:
\begin{equation}\label{Eq:SP1}
x_i \mapsto (-1)^f(-a)^{\bar a_i'}x,
\qquad i=1,\dots,2k~~,
\end{equation}
to write one-variable series  $\bar P_{\bar C_{K_f}}(x)$. However, we {\it propose} the following substitution for the twist knot $K_{p_1<0} \equiv L(2p_1+1,1,1)$, $K_{p_1>0} \equiv L(2p_1-1,1,1)$ 
and double twist knots $K(p_1+1,p2+1)\equiv L(2p_1+1,2p_2+1,1)$:
\begin{equation}
\begin{split}
x_i \mapsto\;& (-1)^f(- \mathrm {sgn} (p_1)a)^{-\,\mathrm{sgn}(p_1)\bar a_i}
a^{2(|p_1|+\mathrm{sgn}(p_1))+2(p_2+1)+1}  \\
&\times q^{q_i+2p_2}\,x,
\qquad i=1,\dots,2k
\end{split}
\label{Eq:SP2}
\end{equation}
where $\mathrm{sgn}(p_1)$ denotes the sign of $p_1$. 
The object that will eventually admit a path-counting interpretation is not the specialized series itself, but the following version:
\begin{equation}
\begin{split}
y(x,a,q)
&=
\frac{\bar P_{\bar C_{K_f}}(qx)}
     {\bar P_{\bar C_{K_f}}(q^{-1}x)}
=
\sum_{l\geq 0}N_l(a,q)x^l \\
&=
\sum_{i,j,l} n_{i,j,l}\,a^i q^j x^l .
\end{split}
\label{Eq:path_counting_GS}
\end{equation}
Under the substitutions (\ref{Eq:SP1},\ref{Eq:SP2}), we have shown that the coefficients $n_{i,j,l}$ are indeed positive integers for some twist knots and double twist knots. This is the key feature that makes a combinatorial interpretation plausible. Our goal is to construct a lattice-path model in which the path counting reproduces these integers.

Clearly, the above exercise on combinatorial interpretation is plausible only if we know the quivers for the 3-pretzel family. Such quivers have been elegantly achieved where full-twist operation on a knot is mapped to unlinking/linking operations on the corresponding augmented quiver \cite{Chauhan:2025iwt}. We will present the necessary details of this procedure in the following subsection.
\subsection{Full twists on knots, unlinking/linking operation on quivers:}
 For a knot $K_1$ containing the twist region $\tau$, adding full twists in $\tau$ produces tower of knots
$\{K_i\}_{i\geq 1}^{\tau},$
where $K_{i+1}$ is obtained from $K_1$ by adding $i$ full twists in $\tau$ (see Figure. 2 of \cite{Chauhan:2025iwt} for details). The family of twist knots $K_{p<0}$ is a standard example, starting with $K_1=4_1$, repeated addition of full twists gives
$K_2=6_1,\, K_3=8_1,\, \ldots.$

For knots related by full twists in a fixed twist region, the $r$-colored HOMFLY-PT polynomials satisfy the recursion relations :
$$
   \sum_{i=0}^{r+1}(-1)^i\,\epsilon_i(\lambda_0,\lambda_1,\ldots,\lambda_r)\,
P_r(K_{j+r+1-i})=0,\qquad j\geq 0, 
$$
where $\lambda_0, \lambda_1,\lambda_2, \ldots$ are the braiding eigenvalues for full-twists between the two strands in the twist region. Here $\epsilon_i$ denotes the $i$-th elementary symmetric polynomial  (see \cite{Chauhan:2025iwt} for details) involving $\lambda_i$'s.
Thus the recursive relation gives the $P_r(K_i)$ for the tower of knots. In fact, the twisting operation on the knot side corresponds on the quiver side as local operations called {\it unlinking} and {\it linking} \cite{Chauhan:2025iwt}.

The unlinking operation enlarges a quiver by one extra node. If $Q=(C,\mathbf{x})$ is a symmetric quiver with $m$ nodes  then unlinking the pair of nodes $(i,j)$ produces a new quiver
\begin{equation}\label{Eq:Unlinking}
Q'(C', \mathbf x, x_{m+1})=U{(i,j)}Q(C,\mathbf x)
\end{equation}
with one extra node. The new adjacency matrix $C'$ elements are related to $C$ as follows:
\begin{align}
C'_{ij} &= C_{ij}-1, \nonumber\\
C'_{s,m+1} &= C_{si}+C_{sj}, \quad s\neq i,j, \nonumber\\
C'_{i,m+1} &= C_{ii}+C_{ij}-1, \nonumber\\
C'_{j,m+1} &= C_{ij}+C_{jj}-1, \nonumber\\
C'_{m+1,m+1} &= C_{ii}+C_{jj}+2C_{ij}-1.
\end{align}
The unlinking operation introduces a new variable corresponding to the added extra node. The generator corresponding to the extra node  is:
\begin{equation}
    x_{m+1}=q^{-1}x_ix_j.
\end{equation}
Similarly, linking produces a new quiver 
\begin{equation}
Q'(C', \mathbf x, x_{m+1})=L{(i,j)}Q(C,\mathbf x)
\end{equation}
with one extra node, with the following relations between the adjacency matrix elements:
\begin{align}
C'_{ij} &= C_{ij}+1, \nonumber\\
C'_{s,m+1} &= C_{si}+C_{sj}, \quad s\neq i,j, \nonumber\\
C'_{i,m+1} &= C_{ii}+C_{ij}, \nonumber\\
C'_{j,m+1} &= C_{ij}+C_{jj}, \nonumber\\
C'_{m+1,m+1} &= C_{ii}+C_{jj}+2C_{ij}.
\end{align}
together with the generator corresponding to the extra node is:
\begin{equation}
   x_{m+1}=x_ix_j. 
\end{equation}
A sequence of unlinkng and linking operations are performed on the nodes of the object, known as the augmented quiver
$
\tilde{Q}^{+}_{K_1},
$
to obtain new augmented quivers corresponding to the tower of knots $\{K_i\}$ obtained through twisting operations. The form of the adjacency matrix associated with the augmented quiver $\tilde{Q}^+_{K_1}(C^+,\mathbf x,x_0)$ of knot $K_1$ will include top row and first column to adjacency matrix of quiver $Q_{K_1}(C,\mathbf x)$ with suitable entries as given in Ref. \cite{Chauhan:2025iwt}. Note that, we refer the top row/first column on the augmented quiver by index `0'  with additional generator $x_0$. Hence, unlinking $(i,j)$ could involve $i,j=0,1,\ldots$.
 
We will utilize the quiver generating procedure \cite{Chauhan:2025iwt} for twist knots and double-twist knots. With this data, we focus on the study of lattice counting for some examples in the following section. In particular, we briefly review the  known path model  for
$T(2,2p+1)=3_1,5_1,7_1,\ldots
$ \cite{Panfil:2018sis, Stosic:2024vbd} and then elaborate our lattice path model results on the twist knots and double twist knots.

\section{Examples}\label{Sec:Examples}
\subsection{ \(T(2,2p+1)\) torus knots}
The lattice-path model for the  $T(2,2p+1)$ torus knots is  known in the literature, and we quote the results for completeness (see \cite{Stosic:2024vbd} for details). 

For negative $T(2,2p+1)$ torus knots (mirror image) with framing $2(2p+1)$, the coefficients $n_{i,j,k}$ appearing in the generating series in (\ref{Eq:path_counting_GS}) become positive integers for the following two choices of the polynomial variables: (i)  $\, a=0,\ q=1$ and (ii) $a=1,\ q=1$. These integers have an interpretation in terms of lattice paths in the positive quadrant of the $(x,y)$-plane. More precisely, they count paths starting from origin: $(0,0)$ and ending at $\left ((2p+1)k,\,2k\right)$, subject to the condition that the path remains below the line
$
y=\frac{2}{2p+1}x.
$

\begin{figure*}[t]
\centering
\subfigure[Lattice path without diagonal steps.\label{Fig:Torus_path1}]{
    \includegraphics[width=0.47\textwidth]{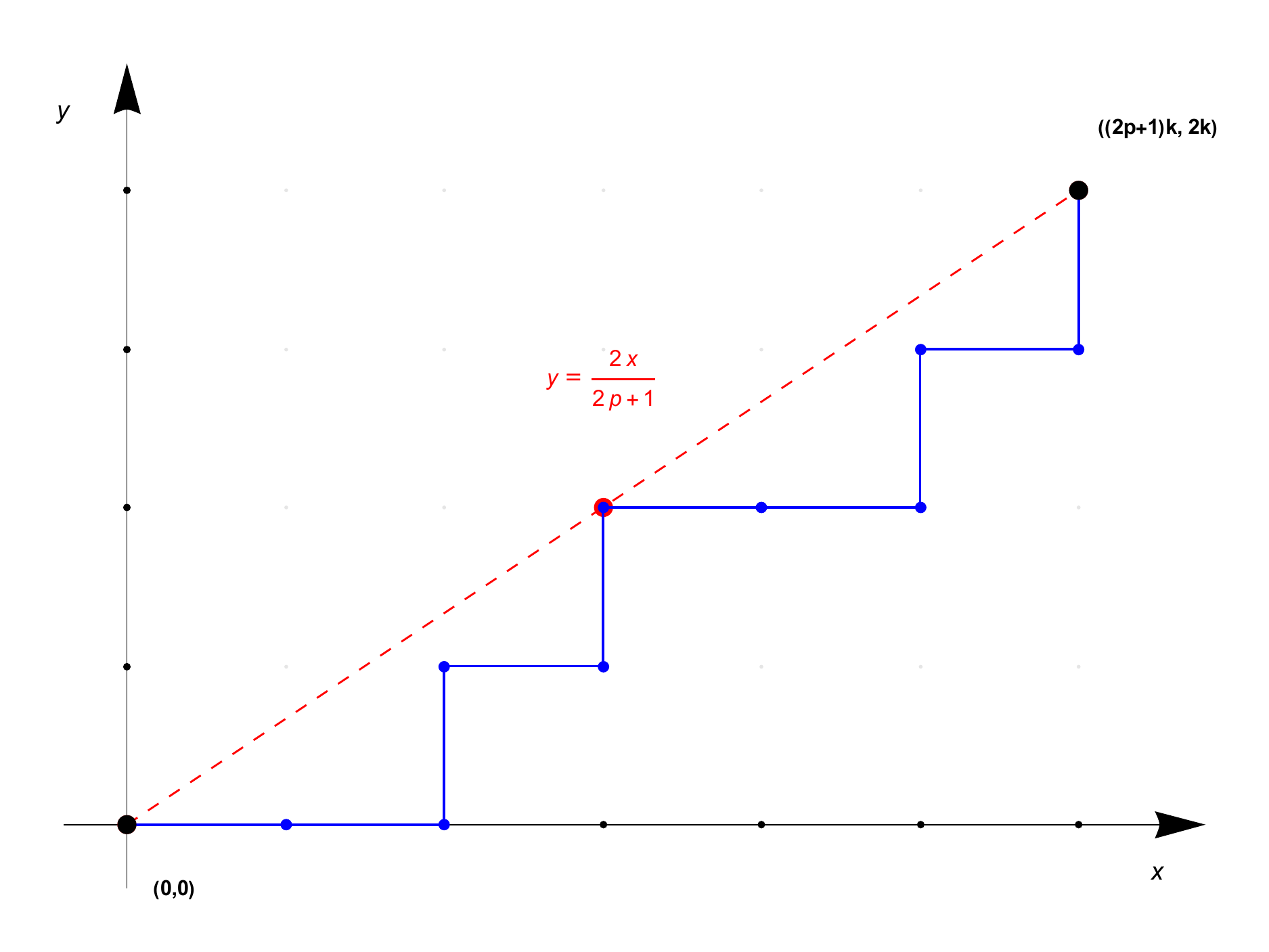}
}
\hfill
\subfigure[Lattice path including the diagonal steps.\label{Fig:Torus_path2}]{
    \includegraphics[width=0.47\textwidth]{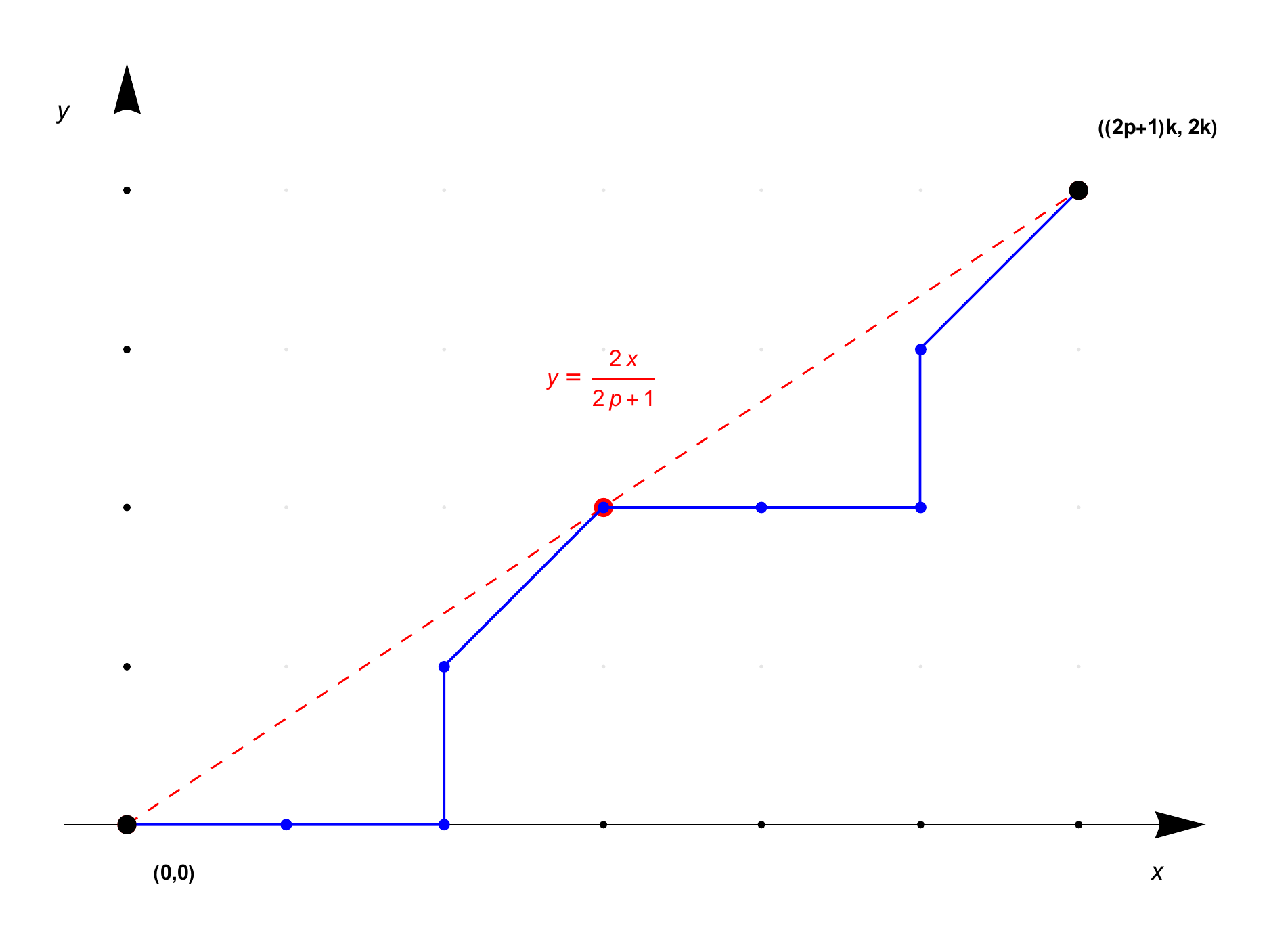}
}
\caption{Lattice paths for $T(2,2p+1)$ torus knots.}
\label{Fig:Torus_path}
\end{figure*}

The choice $a=0$ and $q=1$ leads to the lattice path model with allowed steps $E=(1,0)$ and $N=(0,1)$ \cite{bizley1954derivation}. The corresponding path is shown in Figure \ref{Fig:Torus_path}(\subref{Fig:Torus_path1}). Moreover, the substitution $a=1$ and $q=1$ requires an additional diagonal step $D=(1,1)$ besides the $E=(1,0)$ and $N=(0,1)$ steps. Such a path-counting model is illustrated in Figure ~\ref{Fig:Torus_path}(\subref{Fig:Torus_path2}).
We would like to emphasize the dependence on the twist parameter $p$  in this lattice path counting for $T(2,2p+1)$ torus knot family. In fact, $p$ appears in both the slope of the bounding line as well as the endpoint of the paths. This observation motivates us to look for lattice-path models for the twist knots and double twist knots, and to concisely write twist parameter(s) dependent paths.
In the following subsections, we elaborate our lattice path counting results on the twist knot and double twist knots.

\subsection{$K_{p<0}\equiv L(2p+1,1,1)$:  $4_1,6_1,8_1,\ldots$}

We start with the figure-eight knot, i.e., $K=4_1$. 
In fact, we can generate the tower of knots $$K_{p=-2,-3,\ldots}=6_1,8_1,\ldots$$ by the successive addition of full-twists in the twist region of knot $4_1$. 

\noindent
The adjacency matrix for the augmented quiver $\tilde{Q}^+_{4_1}(C^+_{4_1},\mathbf x\equiv x_0,x_1,\ldots x_5)$ and the vectors $\mathbf a$ and $\mathbf q$ \cite{Chauhan:2025iwt}:
\begin{equation}
\begin{split}
C^{+}_{4_1}
&=
\left[
\begin{array}{cccccc}
1 & 0 & 1 & 1 & 1 & 1 \\
0 & 0 & -1 & -1 & 0 & 0 \\
1 & -1 & -2 & -2 & -1 & 0 \\
1 & -1 & -2 & -1 & -1 & 0 \\
1 & 0 & -1 & -1 & 1 & 1 \\
1 & 0 & 0 & 0 & 1 & 2
\end{array}
\right], \\[2mm]
\bm{a}
&=(a_0,a_1,\ldots,a_5)=(2,0,-2,0,0,2), \\[1mm]
\bm{q}
&=(q_0,q_1,\ldots,q_5)=(-1,0,0,-2,2,0).
\end{split}
\label{Aug:4_1}
\end{equation}
includes the subquiver corresponding to $\tilde{Q}_{4_1}$ when we remove 
the first row and first column.
 The augmented quiver for other knots of this twist family are obtained by unlinking operations (\ref{Eq: Unlinking}):
\begin{equation}
\begin{split}
\tilde{Q}_{K_i}^{+}
={}&
U(0,3i-1)\,U(0,3i-2)\,U(0,3i-3)\,
\\
&U(0,3i-4)\,\tilde{Q}_{K_{i-1}}^{+}, i=2,3,\ldots.
\end{split}
\label{Eq: Unlinking}
\end{equation}
Removing the top row and first column from (\ref{Eq: Unlinking}), we obtain the exact quiver matrix required to begin our calculation. 
\subsubsection*{$4_1$ knot}
We follow Section \ref{Sec:Math_prelm} and obtain the generating series in eqn.~(\ref{Eq:path_counting_GS}) for various framings $f$ of the  $4_1$ knot with negative orientation\footnote{$4_1$ is an amphichiral knot and taking the negative orientation does not change the results. However, in the whole article, we chose the negative oriented knots for uniformity.}. Substituting the generators $x_i$'s  by our proposed specialization (\ref{Eq:SP2}), we find the series $y(x,a,q)$ (\ref{Eq:path_counting_GS}) for $a=0,q=1$. We present the lattice path counting results for some values of framing $f$ in the Table \ref{Nkaq-f:4_1} where the coefficients $N_k(a=0,q=1)$  are naturally identified with a Fuss--Catalan family \cite{aval2007multivariatefusscatalannumbers, oeis}.

 \begin{table}[h!]  
\centering
\begin{tabular}{|c|c|c|}
\hline
\hline
 Framing:$f$ & Sequence&$N_k(a=0,q=1)$ \\
 \hline

$0$& $1, 1, 3, 12, 55, \ldots $&$ \frac{1}{(2k+1)}\binom{3k} {k}$ \\
\hline
 $1$ &$1, 1, 2, 5, 14, \ldots $&$\frac{1}{(k + 1)}\binom{2k}{k}$ \\
\hline
\end{tabular}
\caption{$N_k(a=0,q=1)$ for $4_1$ knot with  framing $f$.}
\label{Nkaq-f:4_1}
\end{table}

For framing $f=0$, the sequence $1,1,3,12,55,\ldots$
is counted by the Fuss-Catalan numbers
\begin{equation}
    N_k^{(0)}(a=0,q=1)
=\frac{1}{2k+1}\binom{3k}{k}.
\end{equation}
Equivalently, $N_k^{(0)}$ counts lattice paths from $(0,0)$ to $(2k,k)$ using steps $E=(1,0)$ and  $N=(0,1)$ which remain weakly below the line $y=\frac{x}{2}.$

Furthermore, for framing $f=1$, the sequence 1,1,2,5,14,\ldots is counted by the ordinary Catalan numbers 
\begin{equation}
    N_k^{(1)}(a=0,q=1)
=\frac{1}{k+1}\binom{2k}{k}.
\end{equation}
Here, $N_k^{(1)}$ is the lattice path counting from $(0,0)$ to $(k,k)$ using the steps $E=(1,0)$ and $N=(0,1)$,
which remain weakly below the line $y=x.$

\subsubsection*{$6_1$ knot}
The augmented quiver for $6_1$ knot is obtained by choosing $i=2$  and applying the successive unlinking operations(\ref{Eq: Unlinking}). That is,
$$\tilde{Q}_{6_1}^{+}
=
U(0,5)\,U(0,4)\,U(0,3)\,U(0,2)\,Q_{4_1}^{+},
$$

We follow the similar procedure as the $4_1$ knot case to obtain  $y(x,a=0,q=1)$ (\ref{Eq:path_counting_GS}). The sequence $N_k(a=0,q=1)$  for the various values of framing $f$ is given in  Table \ref{Nkaq-f:6_1}.

\begin{table}[h!]
\centering
\begin{tabular}{|c|c|c|}
\hline
 Framing:$f$ &Sequence &$N_k(a=0,q=1)$  \\
 \hline
 $0$& $1, 3, 18, 136, 1155, \ldots $&$ 
\frac{3}{5k+3}\binom{5k+3}{k}
$ \\
\hline
 $1$ &$1, 3, 15, 91, 612, \ldots $&$
\frac{3}{4k+3}\binom{4k+3}{k}
$ \\
\hline
  $2$ &$1, 3, 12, 55, 273, \ldots $&$
\frac{3}{3k+3}\binom{3k+3}{k}
$
 \\
\hline
$3$ &$1, 3, 9, 28, 90, \ldots $&$
\frac{3}{2k+3}\binom{2k+3}{k}
$
 \\
\hline
 $4$&$1, 3, 6, 10, 15, \ldots $ &$
\frac{3}{k+3}\binom{k+3}{k}
$ 
\\
\hline
\end{tabular}
\caption{$N_k(a=0,q=1)$  for $6_1$ knot with framing $f$}
\label{Nkaq-f:6_1}
\end{table}
For the above framing values, we can write a concise closed form $N_k^{(f)}$ as:
\begin{equation}
    N_k^{(f)}(a=0,q=1)=\frac{3}{(5-f)k+3}\binom{(5-f)k+3}{k},
\end{equation} 
which counts lattice paths from $(0,0)$ to $((5-f)k+3,k)$ using steps $E=(1,0)$ and $N=(0,1)$, such that, the admissible paths remain  below the  line $y=\frac{x}{5-f}$.
\subsubsection*{$8_1$ Knot}
We can deduce the augmented quiver (\ref{Eq: Unlinking}) from $\tilde Q_{6_1}^{+}$:
$$
\tilde{Q}_{8_1}^{+}
=
U(0,9)\,U(0,8)\,U(0,7)\,U(0,6)\,\tilde{Q}_{6_1}^{+}.
$$
Using the quiver generating series, we obtain the integer sequence (\ref{Eq:path_counting_GS}) for some values of framings $f$ (See Table \ref{Nkaq-f:8_1}). 
\begin{table}[h!]
\centering
\begin{tabular}{|c|c|c|}
\hline
 Framing:$f$ & Sequence&$N_k(a=0,q=1)$  \\
 \hline
 $0$& $1, 5, 45, 500, 6200, \ldots $&$ 
\frac{5}{7k+5}\binom{7k+5}{k}
$ \\
\hline
 $1$ &$1, 5, 40, 385, 4095, \ldots $&$
\frac{5}{6k+5}\binom{6k+5}{k}
$ \\
\hline
  $2$ &$1, 5, 35, 285, 2530, \ldots $&$
\frac{5}{5k+5}\binom{5k+5}{k}
$
 \\
\hline
$3$ &$1, 5, 30, 200, 1425, \ldots $&$
\frac{5}{4k+5}\binom{4k+5}{k}
$
 \\
\hline
 $4$&$1, 5, 25, 130, 700, \ldots $ &$
\frac{5}{3k+5}\binom{3k+5}{k}
$ 
\\
\hline
$5$&$1, 5, 20, 75, 275, \ldots $ &$
\frac{5}{2k+5}\binom{2k+5}{k}
$ \\
\hline
\end{tabular}
\caption{ $N_k(a,q)$  for     $8_1$ knot with framing $f$} .
\label{Nkaq-f:8_1}
\end{table}
For the $8_1$ knot, the framing-dependent closed form is
\begin{equation}
    N_k^{(f)}(a=0,q=1)=\frac{5}{(7-f)k+5}\binom{(7-f)k+5}{k},
\end{equation}
where $f=0,1,2,3,4,5$.

For these framing numbers, $N_k^{(f)}(a=0,q=1)$ counts lattice paths from $(0,0)$ to $((6-f)k+4,k)$ using steps $E=(1,0)$ and $N=(0,1)$, such that every point of the path satisfies $y\leq \frac{x}{6-f}$.

\par
With the detailed tabulation of the integer sequences dependent on framing $f$ for the twist knots $K_{p=-1,-2,-3}$, we attempted to write a concise form of the lattice path counting $N_k^{(p,f)} (a=0,q=1)$. Interestingly, for $f=0$ we get
\begin{equation}
\begin{split}
N_k^{(p,0)}(a=0,q=1)
={}&
\frac{2|p|-1}{(2|p|+1)k+2|p|-1} \\
&\times
\binom{(2|p|+1)k+2|p|-1}{k}.
\end{split}
\label{Eq:pnegf0ClosedForm}
\end{equation}
 In fact,  the integer
$N_k^{(p,0)}(a=0,q=1)$ counts lattice paths from $(0,0)$ to $(2|p|k+2|p|-2,k)$ using steps $E=(1,0)$ and $N=(0,1)$, such that the path remains weakly below the  line $y=\frac{x}{2|p|}$. 

We also succeeded in writing $f=1$ integer sequences data for these $4_1,6_1,8_1$ knots:
\begin{equation}\label{Eq:{p<0,f=1}_Closed_form}
    N_k^{(p,1)}(a=0,q=1)=\frac{2|p|-1}{2|p|k+2|p|-1}\binom{2|p|k+2|p|-1}{k}.
\end{equation}
This number counts lattice paths from $(0,0)$ to $((2|p|-1)k+2|p|-2,k)$ using steps $E=(1,0)$ and $N=(0,1)$, which remain weakly below the  line $y=\frac{x}{2|p|-1}$. 
We perform a similar investigation for $K_{p>0}$ in the following subsection.
\subsection{$K_{p>0}\equiv L(2p-1,1,1)$ \text{ with} $p>1$: $5_2,7_2,9_2,\ldots$}
The quiver generating series for all the knots belonging to this family can be obtained by a sequence of unlinking operations on the augmented quiver  $\tilde{Q}^+_{5_2}$ associated to the knot $5_2$ \cite{Chauhan:2025iwt}. More precisely, we perform the operation: 
\begin{equation}
\begin{split}
\tilde{Q}_{K_i}^{+}
={}&
U(0,4i-1)\,U(0,4i-2)\,U(0,4i-3) \\
& U(0,4i-4)\,\tilde{Q}_{K_{i-1}}^{+},
\qquad i=2,3,\ldots .
\end{split}
\end{equation}
The integer sequences can be extracted for these knots also from $y(x,a=0,q=1)$ derived from the quiver generating series. The results are tabulated for knots $5_2,7_2,9_2$. Further, we could obtain the lattice path counting interpretation for some framings.
\subsubsection*{$5_2$ knot}
The data presented in the Table \ref{Nkaq-f:5_2} for $5_2$ knots for some framing numbers   
\begin{table}[h!]
\centering
\begin{tabular}{|c|c|c|}
\hline
 Framing:$f$ &Sequence &$N_k(a=0,q=1)$  \\
 \hline
 $0$& $1, 4, 30, 280, 2925, \ldots $&$ 
\frac{4}{(6k+4)}\binom{6k+4}{k}
$ \\
\hline
 $1$ &$1, 4, 26, 204, 1771, \ldots $&$
\frac{4}{(5k+4)}\binom{5k+4}{k}
$ \\
\hline
  $2$ &$1, 4, 22, 140, 969, \ldots $&$
\frac{4}{(4k + 4)}\binom{4k+4}{k}
$
 \\
\hline
$3$ &$1, 4, 18, 88, 455, \ldots $&$
\frac{4}{3k+4}\binom{3k+4}{k}
$
 \\
\hline
 $4$&$1, 4, 14, 48, 165, \ldots $ &$
\frac{4}{2k+4}\binom{2k+4}{k}
$ 
\\
\hline
\end{tabular}
\caption{$N_k(a,q)$  for  $5_2$ knot with framing $f$.}   
\label{Nkaq-f:5_2}
\end{table}
can be written uniformly as a Raney family $R_{m,n}(k)$ \cite{ beagley2015raney, michaels1991applications}. Namely,
\begin{equation}
\begin{split}
N_k^{(f)}
&=R_{6-f,\,4}(k)
=\frac{4}{(6-f)k+4}\binom{(6-f)k+4}{k}, \\
&\qquad f=0,1,2,3,4,5.
\end{split}
\end{equation}

Indeed, we can give a lattice path interpretation for the above integer sequences.
$N_k^{(f)}$ counts lattice paths from
$
(0,0)\,\text{to}\,((5-f)k+3,k)
$
with steps $E=(1,0)$ and $N=(0,1)$ which stay weakly below the line
$
y=\frac{x}{5-f}.
$

\subsubsection*{$7_2$ knot}
The integer sequences are again belonging to the Raney family for knot $7_1$ as well and are presented in the Table \ref{Nkaq-f:7_2}. 
\begin{table}[h!]
\centering
\begin{tabular}{|c|c|c|}
\hline
 Framing:$f$ &Sequence &$N_k(a=0,q=1)$ \\
 \hline
 $0$& $1, 6, 63, 812, 11655, \ldots $&$ 
 \frac{6}{(8k+6)}\binom{8k+6}{k}
$ \\
\hline
 $1$ &$1, 6, 57, 650, 8184, \ldots $&$
\frac{6}{(7k+6)}\binom{7k+6}{k}
$ \\
\hline
  $2$ &$1, 6, 51, 506, 5481, \ldots $&$
 \frac{6}{(6k+6)}\binom{6k+6}{k}
$
 \\
\hline
$3$ &$1, 6, 45, 380, 3450, \ldots $&$
\frac{6}{(5k + 6)}\binom{5k + 6}{k}
$
 \\
\hline
 $4$&$1, 6, 39, 272, 1995, \ldots $ &$
\frac{6}{4k+6}\binom{4k + 6}{k}
$ 
\\
\hline
  $5$ &
 $1, 6, 33, 182, 1020, \ldots $ &$
\frac{6}{3k+6}\binom{3k+6}{k}
$ \\
\hline
\end{tabular}
\caption{$N_k(a,q)$  for  $7_2$ knot with framing $f$.}  .
\label{Nkaq-f:7_2}
\end{table}
The general compact form for these framing numbers is
\begin{equation}
\begin{split}
N_k^{(f)}
&=R_{8-f,\,6}(k)
=\frac{6}{(8-f)k+6}\binom{(8-f)k+6}{k}, \\
&\hspace{1.8cm} f=0,1,2,3,4,5.
\end{split}
\end{equation}

For framing $f$, the number $N_k^{(f)}$ counts lattice paths from
$
(0,0)\quad\text{to}\,((7-f)k+5,k)
$
with steps $E=(1,0)$ and $N=(0,1)$ which stay weakly below the line
$
y=\frac{x}{7-f}.
$
\subsubsection*{$9_2$ knot}
Table \ref{Nkaq-f:9_2} gives the integer sequences for $9_2$ knot with framing $f$ 
which also has a neat, concise $f$ dependent closed form:
\begin{table}[h!]
\centering
\begin{tabular}{|c|c|c|}
\hline
 Framing:$f$ &Sequence &$N_k(a=0,q=1)$ \\
 \hline
 $0$& $1, 8, 108, 1776, 32430, \ldots $&$ 
 \frac{8}{(10k+8)}\binom{10k+8}{k}
$ \\
\hline
 $1$ &$1, 8, 100, 1496, 24682, \ldots $&$
\frac{8}{(9k+8)}\binom{9k+8}{k}
$ \\
\hline
  $2$ &$1, 8, 92, 1240, 18278, \ldots $&$
 \frac{8}{(8k+8)}\binom{8k+8}{n}
$
 \\
\hline
$3$ &$1, 8, 84, 1008, 13090, \ldots $&$
\frac{8}{(7k+8)}\binom{7k+8}{k}
$
 \\
\hline
 $4$&$1, 8, 76, 800, 8990, \ldots $ &$
\frac{8}{(6k+8)}\binom{6k+8}{k}
$ 
\\
\hline
  $5$ &
 $1, 8, 68, 616, 5850, \ldots $ &$
\frac{8}{(5k + 8)}\binom{5k + 8}{k}
$ \\
\hline
\end{tabular}
\caption{$N_k(a,q)$  for  $9_2$ knot with framing $f$.}   .
\label{Nkaq-f:9_2}
\end{table}
\begin{equation}
\begin{split}
N_k^{(f)}
&=R_{10-f,\,8}(k)
=\frac{8}{(10-f)k+8}\binom{(10-f)k+8}{k}, \\
& f=0,1,2,3,4,5.
\end{split}
\end{equation}
These integers count the lattice paths from
$
(0,0)\,\text{to}\,((9-f)k+7,k)
$
with steps $E=(1,0)$ and $N=(0,1)$ which stay weakly below the line
$
y=\frac{x}{9-f}.
$

With these detailed computations for positive twist knots $K_{p=2,3,4}$, we are in a position to write the $p$-dependent integer sequence and lattice counting for some framing $f$. In fact, for $f=1$ we deduce the following
$p$-dependent Raney family 

\begin{equation}
\begin{split}
N_k^{(p,f=1)}
&=R_{2p+1,\,2p}(k) \\
&=\frac{2p}{(2p+1)k+2p}
\binom{(2p+1)k+2p}{k}.
\end{split}
\label{Eq:pposf1ClosedForm}
\end{equation}
In this notation, $N_k^{(p,f=1)}$ counts lattice paths from $(0,0)
\, \text{to} \,
(2pk+2p-1,k)
$
with steps $E=(1,0),\, N=(0,1)$,
which stay weakly below the line
$y=\frac{x}{2p}.$

Can we extend our detailed lattice path counting description on twist knots to the double twist knot family?
We attempt a class of double twist knots $K(p+1,2)$ which belongs to the 3-pretzel family $L(2p+1,3,1)$ in the following subsection.
\subsection{ Double twist knots $K(p+1,2)\equiv L(2p+1,3,1)$: with $p\geq 1$: $7_4, 9_5,\ldots$} 
For $p=1$, the double twist knot is $7_4$ whose 3-pretzel notation is $L(3,3,1)$.
Hence, starting from $\tilde{Q}^{+}_{K_1}\equiv \tilde{Q}^{+}_{7_4}$ whose augmented quiver is given in Ref. \cite{Chauhan:2025iwt}, we  obtain the augmented quiver for  family
\(L(2p+1,3,1)|_{p\geq 2}=9_5,\ldots\) using the following unlinking operations: 
$$
\tilde{Q}^{+}_{K_{i+1}}
=
U(0,16i)\,U(0,16i-1)\cdots U(0,16i-7)\,
\tilde{Q}^{+}_{K_i},
\, i\ge 1 .
$$
 The coefficients $N_k(a=0,q=1)$ of the generating series  $y(x,a=0,q=1)$ (\ref{Eq:path_counting_GS}) give rise to the lattice path counting model for both $7_4$ and $9_5$ knots, which we will present in the following subsections.
\subsection*{$7_4$ knot}
The integer sequence for the framings $f=0,1,2,3,4,5$, in case of $7_4$ knot, is presented in Table \ref{Nkaq-f:7_4}. This again corresponds to the Raney family $R_{8-f,4}(k)$ with the following explicit closed form: 
 \begin{equation}
     N_k^{p=1,f}(a=0,q=1)=\frac{4}{(8-f)k+4}\binom{(8-f)k+4}{k}
 \end{equation}
 
 \begin{table}[h!]
\centering
\begin{tabular}{|c|c|c|}
\hline
 Framing:$f$ &Sequence &$N_k(a=0,q=1)$ \\
 \hline
 $0$& $1, 4, 38, 468, 6545, \ldots $&$ 
 \frac{4}{(8k+4)}\binom{8k+4}{k}
$ \\
\hline
 $1$ &$1, 4, 34, 368, 4495, \ldots $&$
\frac{4}{(7k+4)}\binom{7k+4}{k}
$ \\
\hline
$2$& $1, 4, 30, 280, 2925, \ldots $&$ 
\frac{4}{(6k+4)}\binom{6k+4}{k}
$ \\
\hline
 $3$ &$1, 4, 26, 204, 1771, \ldots $&$
\frac{4}{(5k+4)}\binom{5k+4}{k}
$ \\
\hline
  $4$ &$1, 4, 22, 140, 969, \ldots $&$
\frac{4}{(4k + 4)}\binom{4k+4}{k}
$
 \\
\hline
$5$ &$1, 4, 18, 88, 455, \ldots $&$
\frac{4}{3k+4}\binom{3k+4}{k}
$
 \\
\hline
\end{tabular}
\caption{$N_k(a,q)$  for  $7_4$ knot with framing $f$.}   .
\label{Nkaq-f:7_4}
\end{table}
Interestingly, $N_k^{p=1,f}(a=0,q=1)$ counts the  lattice paths from $(0,0)$ to $((7-f)k+3,k)$ using the steps $E=(1,0)$ and $N=(0,1)$ and staying weakly below the  line $y=\frac{x}{7-f}$.

\subsubsection*{$9_5$ knot}

For the values of framing $f=0,1,2,3,4,5$,  Table \ref{Nkaq-f:9_5} contains the integer sequence generated by (\ref{Eq:path_counting_GS}). We observe that these sequences have the framing $f$ dependent form:
\begin{equation}
    N_k^{p=2,f}(a=0,q=1)=\frac{6}{(10-f)k+6}\binom{(10-f)k+6}{k}.
\end{equation}

Here, $N_k^{p=2,f}(a=0,q=1)$ counts the lattice paths from $(0,0)$ to
$((9-f)k+5,k)$ using steps $E=(1,0)$ and $N=(0,1)$, constrained to remain
weakly below the  line $y=\frac{x}{9-f}$. 

\begin{table}[h!]
\centering
\begin{tabular}{|c|c|c|}
\hline
 Framing:$f$ &Sequence &$N_k(a=0,q=1)$ \\
 \hline
 $0$& $1, 6, 75, 1190, 21285, \ldots $&$ 
 \frac{6}{(10k+6)}\binom{10k+6}{k}
$ \\
\hline
 $1$ &$1, 6, 69, 992, 15990, \ldots $&$
\frac{6}{(9k+6)}\binom{9k+6}{k}
$ \\
\hline
  $2$ &$1, 6, 63, 812, 11655, \ldots $&$ 
 \frac{6}{(8k+6)}\binom{8k+6}{k}
$
 \\
\hline
$3$ &$1, 6, 57, 650, 8184, \ldots $&$
\frac{6}{(7k + 6)}\binom{7k + 6}{k}
$
 \\
\hline
\end{tabular}
\caption{$N_k(a,q)$  for  $9_5$ knot with framing $f$.} .
\label{Nkaq-f:9_5}
\end{table}

Furthermore, in the  framing $f=0$, we can combine the  sequences $N_k^{(p=1,2)}(a=0,q=1)$ in a $p$ dependent form as: 
\begin{equation}
\begin{split}
N_k^{(p,0)}(a=0,q=1)
&=R_{2p+6,\,2p+2}(k) \\
&=\frac{2p+2}{(2p+6)k+2p+2} \\
&\quad \times \binom{(2p+6)k+2p+2}{k}.
\end{split}
\label{Eq:DoubleTwistf0ClosedForm}
\end{equation}
    
The integers $N_k^{(p,0)}(a=0,q=1)$  count lattice paths from $(0,0)$ to $((2p+5)k+2p+1,k)$ using steps $E=(1,0)$ and $N=(0,1)$, such that the path remains weakly below the line $y=\frac{x}{2p+5}$.

In this section, we have obtained the path counting model for some twist knots and double twist knots with framing $f$ only when the variables $a=0$ and $q=1$. Our hunt to deduce such lattice paths when $a=1,q=1$ has not been plausible. We hope to succeed in finding such paths in the future.
\section{Summary and Conclusion}\label{Sec:conclusion}
The quiver generating series available for a class of double twist knots \cite{Chauhan:2025iwt} belonging to the $3$-pretzel knot family $L(t_1,t_2,t_3)$ inspired us to investigate lattice path interpretations. Under suitable specializations of the quiver generators $x_i$  (\ref{Eq:SP2}) and after imposing $a=0$ and $q=1$, the corresponding generating series $y(x, a,q)$ (\ref{Eq:path_counting_GS}) gives rise to explicit lattice path counting models for some framings $f$.
\begin{enumerate}
   \item For the negative twist knot family $K_{p<0}$ with framing $f$ we obtained the following results:
\begin{itemize}
    \item $N_k^{(p,f=0)}(a=0,q=1)$ (\ref{Eq:pnegf0ClosedForm}) which counts lattice paths from $(0,0)$ to $(2|p|k+2|p|-2,k)$ with steps $E=(1,0)$ and $N=(0,1)$ lying weakly below the line $y=\frac{x}{2|p|}$.
    
    \item $N_k^{(p,f=1)}(a=0,q=1)$ (\ref{Eq:{p<0,f=1}_Closed_form}) which counts lattice paths from $(0,0)$ to $((2|p|-1)k+2|p|-2,k)$ lying weakly below the line $y=\frac{x}{2|p|-1}$.
\end{itemize}
We observe that the twist parameter $p$ controls the closed-form $N_k^{(p,f)}(a=0,q=1)$ expression as well as the rational slope of the associated lattice-path model. 
\item For positive twist knots $K_{p>1}$ with framing $f=1$, we obtain the following lattice interpretation:
\begin{itemize}
    \item $N_k^{(p,f=1)}(a=0,q=1)$ (\ref{Eq:pposf1ClosedForm}) count lattice paths from $(0,0)$ to $(2pk+2p-1,k)$, again using steps $E=(1,0)$ and $N=(0,1)$, and staying weakly below the line $y=\frac{x}{2p}$.
\item The data of various framing is presented in Tables \ref{Nkaq-f:5_2}, \ref{Nkaq-f:7_2} and \ref{Nkaq-f:9_2}.
    \end{itemize}
 We showed that the resulting sequences fit naturally into families of generalized Fuss-Catalan numbers, or equivalently, Raney numbers. 
\item The integer sequences and lattice path interpretation for the double twist knots $7_4$ and $9_5$, referring to $p=1,2$ respectively in the 3-pretzel knots family $L(2p+1,3,1)$ are :
\begin{itemize}
\item For the fixed framing $f=0$, $N_k^{(p,0)}(a=0,q=1)$ (\ref{Eq:DoubleTwistf0ClosedForm}). This is the Raney number $R_{2p+6,2p+2}(k)$ which counts paths from $(0,0)$ to $((2p+5)k+2p+1,k)$ with steps $E=(1,0)$ and $N=(0,1)$, constrained to lie weakly below the rational-slope line $y=\frac{x}{2p+5}$.  
\item Tables \ref{Nkaq-f:7_4}, \ref{Nkaq-f:9_5} gives integer sequences for various framings. 
\end{itemize}

The results on $T(2,2p+1)$ torus knots when $a=1,q=1$ [see Figure \ref{Fig:Torus_path}(\subref{Fig:Torus_path2})] involve an additional diagonal path in the lattice counting. Do such lattice counting paths derivable for twist knots and double twist knots?
Unfortunately, we are not able to attempt path counting interpretation from the series $y(x,a=1,q=1)$ (\ref{Eq:path_counting_GS}) for the twist knots and double twist knots. We hope to take up this task in the future.

\end{enumerate}
\bibliography{ref}

\end{document}